\input amstex
\documentstyle{amsppt}
\magnification=\magstep0
\define\cc{\Bbb C}

\define\r{\Bbb R}

\define\N{\Bbb N}

\define\jj{\Bbb J}

\define\A{\Cal A}
\define\h{\Cal D}
\define\E{\Cal E}

\define\f{\Cal S}

\define\la{\lambda}
\define\om{\omega}

\define\e{\varepsilon}
\define\va{\varphi }
\define\CB#1{\Cal C_b(#1)}
\define\st{\subset }

\topmatter
 \title
 Spectral characterization of absolutely regular  vector-valued distributions.
  \endtitle
 \endkeywords
 \author
  Bolis Basit and  Alan J. Pryde
\endauthor
 \abstract
{We study  the reduced Beurling spectra $sp_{\Cal {A},V} (F)$  of functions
  $F \in L^1_{loc}
 (\jj,X)$ relative to certain function spaces $\Cal{A}\st L^{\infty}(\jj,X)$ and $V\st L^1 (\r)$, where $\jj$ is $\r_+$ or $\r$ and $X$ is a Banach space. We show that if $F$ is bounded or slowly oscillating  on $\jj$ with $0\not\in sp_{\A,\f} (F)$, where $\A$ is $\{0\}$ or $C_0 (\jj,X)$ for example and $\f=\f(\r)$, then $F$ is ergodic.  This result is new even for $F\in BUC(\jj,X)$ and $\A= C_0(\jj,X)$. If $F$ is ergodic and belongs to  the space $ \f'_{ar}(\jj,X)$ of absolutely regular distributions  and if $sp_{C_0(\jj,X),\f} (F)=\emptyset$, then $\frak{F}*\psi \in C_0(\r,X)$ for all $\psi\in \f(\r)$. Here, $\frak{F}|\jj =F$ and $\frak{F}|(\r\setminus\jj) =0$.
  We show that tauberian theorems for Laplace
transforms follow
   from results about the reduced spectrum.    Our results are more widely applicable than those of previous authors.   We demonstrate this and the sharpness of our results  through  examples}
\endabstract
\endtopmatter
\rightheadtext{  reduced spectra}
\leftheadtext{B. Basit,  A. J. Pryde}
  \TagsOnRight
\document
\pageno=1 \baselineskip=18pt

\head{\S 0. Introduction}\endhead

The goal of this paper is to study the asymptotic behaviour of certain locally integrable functions $F: \jj\to X $ where $\jj$ denotes $\r$ or $\r_+$ and $X$ is  a complex Banach space. The study is motivated by tauberian theorems and their relevance to the behaviour of solutions of Cauchy problems in Banach spaces as in [2], [3], [4, Chapter 5], [6]. When $\jj=\r_+$, the term   ``tauberian"  has been used to describe theorems where  the asymptotic behaviour of a function is deduced  from properties of its Laplace transform, or equivalently its Laplace spectrum $sp^{\Cal{L}} (F)$ (see [21], [4, p. 275]). Improvements were made by employing  a smaller spectrum, the weak Laplace spectrum $sp^{w\Cal{L}} (F)$ (see [4, p. 324], [16], [17]). In this paper we use an even smaller spectrum,
  the reduced (Beurling) spectrum $sp_{\A} (F)$ of $F$ relative
 to   various closed subspaces $\A$ of $L^{\infty}(\jj,X)$. This  spectrum has been used before in this context  (see [5], [6], [7], [16]). It has   the advantage that it unifies the two cases $\jj=\r_+$ and $\jj=\r$.  Moreover, we are able to consider functions whose Fourier transforms  are not regular distributions and  avoid some geometrical restrictions  on  $X$ that were imposed in [17] for example.  Importantly, spectral criteria  for solutions of evolution equations  are readily  related  to  reduced spectra  (see  Theorem 3.5, [6]).
\newline \indent
In section 1 we describe our
notation  and prove some preliminary results.
\newline \indent In section 2 we consider a more general spectrum $sp_{\A,V} (F)$, the reduced spectrum of $F$ relative to $(\A,V)$, where $V\st L^1(\r)$, a spectrum first studied in [9]. Typically, $V$ is one of the spaces $\h=\h(\r)$, $\f=\f(\r)$ or $L^1(\r)$. If  $F\in L^{\infty}(\r,X)$, $\A=  \{0\}$ and $V=L^1(\r)$, then $sp_{\A,V} (F)=sp^{B} (F)$ the classical Beurling spectrum.
    We study the conditions imposed on $\A$ and relate them to previous ones in Proposition 2.1. Then we develop  some basic properties of the reduced spectrum in Propositions 2.2 and 2.3. Our main results are stated in Theorems 2.5 and 2.6. The former deals with ergodicity. We show for example
   that if $F$ is bounded or  slowly oscillating  on $\jj$ with $0\not\in sp_{\A,\f} (F)$, where $\A$ is $\{0\}$ or $C_0 (\jj,X)$, then $F$ is ergodic. Theorem 2.6 deals with functions $F \in\f'_{ab}(\jj,X)$, the space of absolutely regular distributions, with $sp_{C_0(\jj,X),\f} (F)$  countable. It is a generalized tauberian theorem providing spectral conditions under which $F$ has various types of  asymptotic behaviour.  For example (Theorem 2.6 (iv)), if $sp_{C_0(\jj,X),\f} (F)$ is countable and non-empty and $\gamma_{-\om} F$ is ergodic for each  $\om\in sp_{C_0(\jj,X),\f} (F)$, then $(\frak {F}*\psi)|\, \jj$ is asymptotically almost periodic for all $\psi\in \f (\r)$.   Versions of Theorem 2.5 and Theorem 2.6 (i), (ii), (iii), (v) are already known when $\jj=\r_+$ and $sp_{\A,\f} (F)$ is replaced by the larger spectrum $sp^{w\Cal{L}} (F)$ (see [17]). Theorem 2.6 (iv), (vi) seem to be new  for any spectrum.
      Proposition 2.7  states that if $F\in L^1_{loc} (\jj,X)$  with $sp_{C_0(\jj,X),\h} (F)=\emptyset$ and if the convolution $(F*\psi)|\,\jj$ is  uniformly continuous for some $\psi\in \h(\r)$
 then $F*\psi \in C_0(\r,X)$. Chill [17, Proposition 2.1]  obtained this same conclusion under the stronger assumptions that $F\in L^1_{loc}(\r,X)$ and $\widehat{F}\in \f'_{ar}(\r,X)$. In particular, if $F\in L^p (\r,X)$ where $1\le p < \infty$, then $F$ satisfies the assumptions of Proposition 2.7. However, as is well-known, when $p >2$ there are functions $F\in L^p (\r,X)$ for which $\widehat{F}$ is not a regular distribution and so the result of [17] does not apply. Even when $1\le p\le 2$ special geometry on $X$ is required in order that every $F\in L^p (\r,X)$ has a Fourier transform which is regular.

In section 3  we  establish some properties of   weak Laplace, Laplace
and Carleman
spectra which are analogous to those of  Beurling  spectra.   Also, if $F\in \f'_{ar} (\r_+,X)$ and ${\Cal{ A}\supset
C_0 ({\Bbb R}_+,X)}$ then
$sp_{\Cal {A},V} (F)\subset sp^{w\Cal{L}}(F)$ (Proposition 3.2).  As a consequence,
   we  strengthen several theorems about the asymptotic
behaviour of absolutely regular tempered distributions, replacing  Laplace and weak Laplace spectra by   $ sp_{C_0(\r_+,X),\f}(F)$ (see  Remark 3.3). In Theorem 3.5 we obtain a spectral condition satisfied by  bounded mild solutions
of the evolution equation
 $\frac{d u(t)}{dt}= A u(t) +  \phi (t) $, $u(0) \in {X}$, $t\in  {\jj}$,
  where
 $A$ is a closed linear operator   on ${X}$ and $\phi\in L^{\infty} (\jj, {X})$.  This  generalizes earlier results  where it is assumed that $u,\phi \in BUC(\jj,X)$ (see [4, Proposition 5.6.7], [6, Theorem 3.3, Corollary 3.4]).

 \head{\S 1. Notation, Definitions and preliminaries}\endhead

In this paper $\r_+=[0,\infty)$, $\jj\in\{\r_+,\r\}$, $\N= \{1,
2,\cdots \}$, $\cc_+
 =\{\la\in \cc: \text{Re \,}\la
>0\}$ and $\cc_-=\{\la\in \cc: \text{Re \,}\la < 0\}$.
  By $X$ we  denote a   complex Banach space. If $Y$, $Z$ are locally
 convex topological spaces, $L(Y,Z)$  denotes the space of all bounded linear operators from $Y$ to
 $Z$.
   The Schwartz spaces of test functions and rapidly decreasing functions are denoted by $\h(\r)$ and $\f
   (\r)$ respectively.  Then  $\h' (\r,X)= L(\h(\r),X)$ is the  space of $X$-valued distributions and   $\f' (\r,X)= L(\f(\r),X)$ is the  space of $X$-valued tempered distributions
   (see [4, p. 482], [30, p. 149] for $X=\cc$). The space  of absolutely regular distributions is defined by

(1.1) \qquad $\f'_{ar} (\jj,X)=\{ H\in L^1 _{loc}(\jj,X): H\va \in L^1(\jj,X)$ for all $\va\in \f(\r) \}$.

 \noindent   The action of an element $ S\in \h' (\r,X)$ or  $\f' (\r,X)$ on
  $\va\in\h(\r) $ or $\f(\r)$ is denoted  by $<S,\va>$.
 If $F$ is an  $ X$-valued function defined on  $\jj$ and $s\in \jj$
then   $F_s$,  $\Delta _sF$,  $|F|$  stand for the functions defined on
$\jj$ by $F_{s}(t) = F(t+s)$, $\Delta _sF (t)= F_{s}(t)-F (t) $ and $|F| (t)= ||F(t)||$.
Also  $||F|| _{\infty} =
  \text {  sup}_{t\in \jj} ||F(t)||$.
 If $F \in L_{loc}^1 (\jj, X)$ and $h>0$, then $PF$,  $M_hF$  and $\check {F}$ (when $\jj=\r_+$)  denote  the
$indefinite\,\, integral$,  $mollifier$  and $reflection$ of $F$ defined
respectively by $PF(t) = \int_{0}^{t} F (s)\,ds$, $M_h F (t)=
(1/h)\int_0^h F (t+s)\,ds$ for $t\in \jj$ and $\check {F} (t) =F(-t)$ for $t\in \r$.
 For $g\in
   L^1(\r)$ and $F\in
   L^{\infty}(\r,X)$ or $g\in
   L^1(\r,X)$ and $F\in
   L^{\infty}(\r)$ the $Fourier\,\, transform$ $\widehat{g}$ and $convolution$ $F*g$ are
   defined respectively by $\widehat{g} (\om)=\int_{-\infty}^{\infty} \gamma_{-\om} (t)\, g(t)\,
   dt$ and $F*g (t)= \int_{-\infty}^{\infty}F(t-s) g(s)\,
   ds$, where  $\gamma_{\om} (t)= e^{i\, \om t} $ for $\om\in\r$.
    The $Fourier\, transform$ of $ S\in \f' (\r,X)$ is
  the tempered distribution $\widehat{S}$  defined by

  (1.2)\qquad  $<\widehat{S},\va> =<S,\widehat {\va}>$\,\, for all $\va\in
   \f(\r)$.

\noindent Set $\widehat{\h}(\r)=\{\widehat{\va}: \,\,\, \va\in \h
(\r)\}$.
    The $Fourier\, transform$ of $F \in L_{loc}^1 (\r, X)$
   is the distribution $\widehat {F}\in L(\widehat{\h}(\r),X)$
   defined by

  (1.3)\qquad $<\widehat{F},\psi> =<F,\widehat{\psi}>$\,\, \, for all $\psi\in
  \widehat{ \h}(\r)$.

\noindent Throughout the paper all integrals are Lebesgue-Bochner
integrals ([4, pp. 6], [19, p. 318],
   [20, p. 76]). All convolutions  are understood as convolutions of functions defined on $\r$. Given

 (1.4) \qquad     $F\in  W (\jj,X) \in  \{L^1_{loc} (\jj,X), \f'_{ar} (\jj,X), L^{\infty} (\jj,X)\}$,

\noindent     we denote by

 (1.5) \qquad  $\frak{F}$ the function given  by  $\frak{F} |\jj =F$ and $\frak{F} |(\r\setminus \jj) =0$.

\noindent Then $\frak{F}\in W(\r,X)$. In addition, if $g\in L^{\infty}_c (\r)=\{f\in L^{\infty} (\r): f $ has compact support $\}$, then for some constant $t_{g}$

 (1.6) \qquad $\frak{F}*g\in W (\r,X) \cap C(\r,X)$ and  if $\jj=\r_+$,  $\frak{F}*g (t)=0$ for all $t\le t_{g}$.

\noindent  It follows that if  $h > 0$ and  $s_h= (1/h)\chi_{(-h,0)}$, where   $\chi_{(-h,0)}$ is the characteristic function of $(-h,0)$, then

(1.7) \qquad  $\frak{F}*s_h\in W(\r,X)\cap C(\r,X)$, \,\,
 $M_h F= (\frak{F}*s_h) |\jj$ and

\qquad\qquad\,\,\, if $\jj=\r_+$,  $ \frak{F} *s_h (t) =0$ for all $t \le -h$.

 We  use convolutions of functions $F\in W=W(\jj,X)$ and $g \in V=V(\r)\in\{ \h(\r), \f(\r), L^1(\r)\}$, with

(1.8) \qquad $V=\h(\r)$ if $W= L^1_{loc} (\jj,X)$, $V=\f(\r)$ if $W=\f'_{ar}(\jj,X)$ and

\qquad \qquad \,\,\,  $V=L^1(\r)$ if $W=L^{\infty}(\jj,X)$.

 \noindent The following properties of the convolution  are repeatedly used
(see [30, p. 156 (4)], [28, 7.19 Theorem (a), (b), pp. 179-180] when $X=\cc$):

\noindent If $ F\in W(\jj,X)$ and $\va\in V(\r) $  with $W,V$ satisfying (1.8), then

(1.9)\qquad  $\frak{F}*\va\in W(\r,X) \cap C
(\r,X)$.

\noindent Indeed, the cases $W=L^1_{loc}(\jj,X)$ and   $W=L^{\infty}(\jj,X)$ are obvious. If   $F\in \f'_{ar}(\jj,X)$, then
$|F|\in \f'_{ar} (\jj,\cc)$. By [22, Theorem. (b)] there is an integer  $k \in \N$  such that

(1.10.) \qquad $\frac{|F|}{w_k}= f\in L^1(\jj)$, where $w_k (t)= (1+t^2)^k$.

\noindent  Using (1.10), we easily conclude (1.9).

 Moreover, if $\psi \in V(\r)$ or $\psi\in L_c^{\infty}(\r)$, then

(1.11)\qquad  $(\frak{F}*\va)*\psi= (\frak{F}*\psi)*\va$.

 Also we need  the  following analogue of Wiener's
theorem on Fourier series.

\proclaim{Lemma 1.1}  Let $f\in L^1 (\r)$  with $\widehat {f} \not =
0$ on a compact set $ K$. Then  there exists $g\in L^1 (\r)$ such that
$\widehat {g} \cdot \widehat {f} =1$ on $ K$. Moreover, one can choose $g$
such that $\widehat{g}$ has compact support and, if $f\in \f (\r)$,
with $g\in \f (\r)$.
\endproclaim

\demo{Proof} Choose a bounded open set $U$ such that $K\st U$ and
$\widehat {f} \not = 0$ on $\overline {U}$ the closure of $ {U}$. By
[15, Proposition 1.1.5 (b), p. 22], there is $k \in L^1(\r)$ such
that $\widehat {k} \cdot \widehat {f} =1$ on $ \overline {U}$. Now, choose
$\va\in \h(\r)$ such that $\va=1$ on $K$ and supp $\va\st
\overline {U}$.  By [28, Theorem 7.7 (b)] there is  $ \psi\in
\f(\r)$ such that $\widehat{\psi}=\va$. Take $g= k*\psi$. Then
$\widehat {g}$ has compact support and if $f\in \f(\r)$, then
$\widehat{g}\in \h(\r)$ and so $g\in \f(\r)$. $\square$
\enddemo

In the following  proposition  $\psi$ will denote an element  of $ \f
(\r)$ with the  properties:

\qquad $\widehat {\psi}$  has compact support,  $\widehat
{\psi} (0)=1$ and $\psi $ is non-negative.

\noindent An example of such  $\psi$ is given by
$\psi=\widehat{\va}^2$, where $\va(t)=a\, e^{\frac{1}{t^2-1}}$ for
$|t|\le 1$,  $\va=0$ elsewhere on $\r$,  with $a$ some suitable
constant.

\proclaim{Proposition 1.2} (i)
 The sequence $\psi_n(t)= n\, \psi (n\,t)$ is an approximate identity
for the space of uniformly continuous functions $UC(\r,X)$, that
is $\lim _{n\to \infty}||u*\psi_n- u||_{\infty}= 0$  for all $u\in UC(\r,X)$.

(ii) $\lim_{h\searrow 0}||M_h u- u||_{\infty} = 0$  for all $u\in UC(\jj,X)$. In particular if $M_h u\in BUC(\jj,X)$ for all $h > 0$ then $ u\in BUC(\jj,X)$.
\endproclaim

\demo {Proof} (i) Given  $u\in UC(\r,X)$ and $\e >0$ there exists
$k >0$ such that  $||u (t+s)-u(t)|| \le k  |s|+\e$ for all
$t,s\in \r$. In particular $u\in \f'_{ar} (\r,X)$.  Also,
 $u*\psi_n (t)- u(t)= \int_{-\infty}^{\infty} [u (t-\frac{s}{n})-u(t)]
 \psi (s)\,ds$ which gives $||u*\psi_n - u||_{\infty}\le (k/n)\int_{-\infty}^{\infty}
 |s|\psi (s)\,ds+ \e \int_{-\infty}^{\infty}
 \psi (s)\,ds$ and  (i) follows.

 (ii) Since $||M_h u- u||_{\infty}\le \text{\, sup} _{t\in \jj, 0 \le s \le h}||u (t+s)- u(t)||$, part (ii) follows.
  \P
\enddemo
\proclaim{Proposition 1.3} (i) Let $F\in W (\r_+,X)$  and $g\in V(\r)$ with $W, V$ satisfying (1.8). Then $\lim_{t\to -\infty} ||\frak{F}*g (t)||=0$.

(ii) Let $F\in L^{\infty} (\r,X)$ and  $F|\r_+ =0$. Then $(F*f)|\r_+ \in C_0 (\r_+,X)$ for each $f\in L^1(\r)$.
\endproclaim

\demo {Proof} Part (ii) and  the cases $F\in L^1_{loc} (\r_+,X)$ and $F\in L^{\infty} (\r_+,X)$ of part (i) can be shown by simple calculations. If $F\in \f'_{ar} (\r_+,X)$, then from (1.10) ${|F|}/{w_k}= f\in L^1(\r_+)$ for some $w_k (t)= (1+t^2)^k$. Since $\va\in \f(\r)$,
 $||w_k \va||_{\infty}=c_k < \infty $.  It follows that
 $||\frak{F}*\va (t)||= ||\int_0^{\infty} \va (t-s) F(s)\, ds|| \le  \int_0^{\infty} |\va| (t-s) |F|(s)\, ds \le  c_k \int_0^{\infty} \frac{w_k (s)}{w_k (t-s)} f(s)\, ds$. Since $\frac{w_k (s)}{w_k (t-s)}\le 1$ for each $t\le 0, s \ge 0$ and $\lim_{t \to -\infty} \frac{w_k (s)}{w_k (t-s)}$

 \noindent $=0$ for each $s \ge 0$, it follows that $\lim_{t\to -\infty} ||\frak{F}*\va (t)||=0$ by the Lebesgue convergence theorem.
 \P
\enddemo
\head{\S 2 Reduced spectra for regular distributions}\endhead

In this section we introduce  the reduced spectrum $sp_{\A,V} (F)$ of a function $F\in L^1_{loc} (\jj,X)$  relative to $\A, V$, where $\A \st L^{\infty} (\jj,X)$ and $V\st L^1(\r)$. We usually impose the following conditions on $\A$.

(2.1) \qquad  $\A$ is a  closed subspace of $
L^{\infty}(\jj,X)$ and is $BUC$-$invariant$; that is

\qquad\qquad \,\,\, if $\phi\in BUC(\r,X)$ and $\phi|\,\jj\in \A$,
then $\phi_a|\,\jj\in \A$ for each $a\in \r$.

The property of being $BUC$-$invariant$ was first introduced in  [5, $(P.\Lambda)$, Definition 1.3.1] and called the Loomis property $(P.\Lambda)$ for classes $\A \st BUC(\jj,X)$. The notion was extended to classes $\A\st L^1_{loc} (\jj,X)$ in  [7, (1.III$_{ub}$)]. In [9], this property was called  $C_{ub}$-invariance.

We note that if $\jj=\r$, then  $\A$ is $BUC$-invariant if and only if $\A \cap BUC (\r,X)$ is a translation invariant subspace of $BUC (\r,X)$.
If $\jj=\r_+$, then $\A$ is  $BUC$-invariant if and only if $\A \cap BUC (\r_+,X)$ is  a positive invariant subspace of $BUC (\r_+,X)$ (that is $\phi_t \in \A$ for all  $\phi\in \A$, $t\ge 0$) with the additional  property that $u\in \A$ whenever $u\in BUC (\r_+,X)$ and  $ u_t \in \A $ for all $t \ge 0$.
  Such subspaces of $BUC(\r_+,X)$ were subsequently called
 $S$-biinvariant  ([16, (1.1), p. 17], [3,\S 2]).

For $\A$ satisfying  (2.1),  $V\st L^1 (\r)$ and $F\in L^1_{loc} (\jj,X)$,  a point $\om\in \r$ is called $(\A, V)$-$regular$
for $ F$ or $\frak{F}$, if there is $\va\in V$ such that $\widehat {\va}(\om)\not =0$ and
$(\frak{F}*\va)|\,\jj \in \A$. The $reduced$ $Beurling$ $spectrum$ of $F$ or $\frak{F}$
relative to $(\A,V)$ is defined by

(2.2) \qquad $sp_{\A,V} (F) =\{\om\in \r: \om$ is not an
$(\A,V)$-regular point for $F \}=$

\qquad\qquad  \,\, $ \{\om\in \r: \va \in V, (\frak{F}*\va)|\, \jj\in \A
$ implies  $\widehat{\va} (\om)=0 \}=sp_{\A,V} (\frak{F})$,

\noindent provided  the convolution $\frak{F}*\va$ and the restriction $(\frak{F}*\va)|\, \jj$ are defined for all $\va\in V$ (see [10, (1.6)]).
Further, if $H\in L^1_{loc} (\r,X)$   we also define (see [9, Definition 3.1])

(2.2$^{*}$) \qquad $sp_{\A,V} (H)=\{\om\in \r: \va \in V, (H*\va)|\, \jj\in \A
$ implies  $\widehat{\va} (\om)=0 \}$.

\noindent It is clear that  $sp_{\A,V} (F)$ and $sp_{\A,V} (H)$ are  closed subsets of $\r$.
 If $\jj=\r$, then  $F= H|\, \r= \frak{F}$ and  so (2.2) and (2.2$^{*}$) give the same spectrum. If $\jj=\r_+$ and  $F=H|\,\r_+$ we are  interested in comparing   $  sp_{\A,V} (F)$ defined by (2.2) with  $ sp_{\A,V} (H)$ defined by (2.2$^{*}$) (see Proposition 2.2).

For  $ F\in L^{\infty}(\jj,X)$ and $V= L^1(\r)$ we write $sp_{\A} (F)= sp_{\A,L^1(\r)} (F)$.

 If $ F\in W(\jj,X)$ and $V= V(\r)$ with $W,V$ satisfying (1.8), then the convolution $\frak{F}*g$ and the restriction $(\frak{F}*g)|\, \jj$ are defined for all $g\in V(\r)$. So, $sp_{\A,V} (F)$    is well defined.

  This is an
extension of the definitions in  [5, (4.1.1)], [6, (2.9)], [16, Definition 1.14, p. 24]. In those references the conditions on $\A$
 are more restrictive and $ F\in L^{\infty}(\r,X)$. In
particular, if $\A=\{0\}$ and $ F\in L^{\infty}(\r,X)$ then
$sp_0 (F)=sp_{\{0\}} (F)$ is the classical Beurling spectrum $sp^{B} (F)$
[27, p. 183]. If  $ F\in
\f'_{ar}(\r,X)$, then $sp_{0,\f} (F)= $ supp
$\widehat{F} = sp^{\Cal{C}} (F) $ (the Carleman spectrum). Indeed, the first equality is straightforward   and the second is proved in  [26, Proposition 0.5]. If $ F\in L^{\infty}(\jj,X)$, then $\frak{F}*f\in BUC(\r,X)$ for all $f\in L^1(\r)$. It follows that \,\,\,
 $sp_{\A} (F)= sp_{\A\cap BUC(\jj,X)} (F)$.

  Our approach of defining   reduced spectra  via convolutions is widely applicable.  For $ F\in BUC(\jj,X)$ and $\A \st BUC(\jj,X)$, there is  also an operator theoretical approach using $C_0$-semigroups and groups. In [18] it is proved that the two approaches are equivalent for such $F$ and $\A$ (see also [9, Theorem 3.10]). In [24] there is  an  unsuccessful attempt to
extend the operator theoretical approach to $ F\in BC(\jj,X)$ and $\A \st BC(\jj,X)$ (see [25]).

The space $\A_{\frak {g}} = \frak{g}\cdot AP(\r,X)$, where
$\frak{g}(t) =e^{it^2}$ for $t\in\r$,  satisfies (2.1) and $\A_{\frak {g}} \cap BUC(\r,X)=\{0\}$. We conclude   that if $0\not = F\in BC(\r,X)$, then   $sp_{\A_{\frak{g}}} (F)= sp^B (F)\not = \emptyset$. In particular, $sp_{\A_{\frak{g}}} (F)\not = \emptyset$ for each $0\not =F\in \A_{\frak{g}}$. A sufficient condition to have the property  $sp_{\A} (F)=\emptyset$   for each $F\in \A \st L^{\infty} (\jj,X)$ is the following inclusion

  (2.3)\qquad  $(\frak{F}*f)|\, \jj\st \A$ for each $F\in \A$ and $ f\in L^1(\r)$.

Note that if $\A\st BUC(\jj,X)$ satisfies (2.1), then using the properties of Bochner integration we find $\A$ satisfies (2.3). The space $\A_{\frak{g}}$  does not satisfy (2.3).

 Examples of spaces  $\A$ satisfying (2.1), (2.3) include (using $\A(\jj,X)=\A(\r,X)|\jj $)

  $\{0\}$,\,  $C_0 (\jj,X)$, \, $ AP(\r,X)$,\, $LAP_b (\r,X)$,
  $ AA(\r,X)$,\,\, $ EAP(\jj,X) $,

$AAP(\jj,X)= AP(\jj,X)$
$\oplus\, C_0 (\jj,X)$,
 $AAA(\jj,X)= AA(\jj,X)\oplus\, C_0 (\jj,X)$,

\noindent  the spaces consisting  respectively of the zero
function (when $\jj=\r$),  continuous functions vanishing at
infinity, almost periodic, Levitan bounded  almost periodic, almost automorphic  functions ([1], [5], [23]), Eberlein (weakly)  almost periodic ([5, Definition 2.3.1]),
asymptotically almost
periodic functions  (when $\jj=\r_+$)([5, Definitions 2.2.1,
2.3.1, (2.3.2)]) and asymptotically almost
automorphic functions.

 For $\la\in \cc_+$  set

\qquad \qquad $f_{\la}(t) = \cases{ e^{-\la t},
\text{\,\, if\,} t \ge 0}
\\ { 0,\,\,\,\text {  if \,} t<
 0}\endcases$ and $f_{-\la}= -\check {f_{\la}}$.

\noindent Then  $f_{\la}, \check {f_{\la}} \in L^1(\r)$ for all $\la\in \cc\setminus i\r$.
 For $\phi\in L^{\infty} (\r,X)$ and $t\in \r$ define

 \qquad \qquad $\Cal {C} {\phi_t} (\la)= \phi* \check{f_{\la}}(t)=$
$ \cases{ \int_{0}^{\infty}e^{-\la s} \phi (s+t)\,ds,
\text{\,\, if\,\,\,\,} \text{Re\,\,}\la > 0}\\ { -\int_{-\infty}^ 0 e^{-\la s} \phi (s+t)\,ds,  \text {\,\,\, if \,} \text{Re\,\,}\la <
 0}\endcases$

\noindent Obviously, $\phi* \check{f_{\la}}\in BUC(\r,X)$ for all $\la\in \cc\setminus i\r$.   We consider the property

  (2.3$^*$)\qquad $(\frak{F}* \check{f_{\la}})|\jj \in \A$ for each $F\in \A$ and $\la\in \cc\setminus i\r$

\noindent as well as the following

(2.4) \qquad $H\in L^1_{loc}(\r,X)$ and $H| (-\infty,0)$ is bounded.

\proclaim{Proposition 2.1} Let  $\A\st L^{\infty}(\jj,X)$ be a closed subspace.

(i) If  $\A$ is $BUC$-invariant  and $\jj=\r_+$, then $C_0 (\r_+,X) \st \A$.  However, this is not necessarily true if  $\jj=\r$.

(ii)  If $\A$ satisfies (2.3), then
  $\A$  is $BUC$-invariant.

(iii)   $\A$ satisfies (2.3$^*$) if and only if $\A$ satisfies (2.3).
\endproclaim
\demo{Proof} (i) By the $BUC$-invariance  of $\A$,  if $F\in BUC(\r,X)$ and $F$ has
compact support in $(-\infty,0]$, then $F_t|\r_+ \in \A$ for all
$t\in \r$.  It follows that  the space of continuous functions with
compact support $C_c (\r_+,X)\st \A$ and so $C_0 (\r_+,X)\st \A$
(see also the proof of Theorem 2.2.4 in [5, p. 13]). A counter-example
 for the case $\jj=\r$ is $\A= AP(\r,X)$.

(ii) The case $\jj=\r$: Since $L^1(\r)$ is translation invariant, the set $\A* L^1(\r)$ is  translation invariant  too. By Proposition 1.2 (i),
  $\A* L^1(\r)$ is a dense subset of  $\A \cap BUC(\r,X)$.  As  $\A$ is  closed, $\A \cap BUC(\r,X)$ is translation invariant.

The case $\jj=\r_+$: By Proposition 2.2 (ii) below, we conclude that
 $(H_t*f) |\,\r_+ = (H*f_t) |\,\r_+\in \A$ for each $H\in L^{\infty} (\r,X)$ such that $H|\r_+\in \A$ and each $f\in L^1(\r)$ and $t\in \r$. If $H \in BUC(\r,X)$, then again using  the approximate identity of
  Proposition 1.2 (i) we conclude that $H_t  |\,\r_+\in \A $ for each $t\in \r$. This  gives (ii).

(iii) Obviously, (2.3) implies (2.3$^*$). For the converse we begin by showing  that  $ E= span \,\{f_{\la}: \text{ Re\,\,}\la \not =0\}$ is a dense  subspace of  $L^1(\r)$.
    Indeed, if $E$ is not dense in $L^1(\r)$, then by the Hahn-Banach theorem there is
$0\not = \phi\in L^{\infty}(\r) = (L^1(\r))^*$ such that
$\Cal  {C}\phi(\la)= \int_0^{\infty} e^{-\la t}\phi (t)\, dt=0$
if  Re $\, \la >0$
and $\Cal  {C}\phi(\la)= -\int_0^{\infty} e^{\la t}\phi (-t)\, dt=0$ if Re $\, \la  < 0$. This means that the Carleman transform $\Cal  {C}\phi$ is zero on $\cc\setminus i\r$ and implies $sp^{C} (\phi)=\emptyset$ and so $\phi =0$ (see [26, Proposition 0.5 (ii)]). This is a contradiction  showing that $E$ is dense in $L^1(\r)$. Given (2.3$^*$) it follows
  that $(\frak{F}*f)|\, \jj \in \A$ for each $F\in \A$ and $f\in E$. Since $E$ is a dense subspace of $L^1(\r)$ and $\A$ is closed, (2.3) follows.
 $\P$
\enddemo

\proclaim {Proposition 2.2}  Let  $\A \st L^{\infty} (\jj,{X})$ be a closed subspace  satisfying (2.3). Assume  that  $H\in W(\r,X)$  satisfies (2.4) if $\jj=\r_+$ and let $F = H|\, \jj$.

(i)   If $W,V$ satisfy (1.8),
  then $sp_{\A,V} (H)= sp_{\A,V} (F)  $.

(ii)  If  $H\in L^{\infty} (\r,X)$ and  $H|\,\jj\in \A$,  then
  $(H*f)|\,\jj \in \A$ for each $f\in L^1(\r)$.

(iii) If   $H\in L^{\infty} (\r,X)$, then
 $sp_{\A,\f} (F)= sp_{\A,\f} (H)=sp_{\A} (H)=$
  $sp_{\A} (F)$. In particular, $sp_{0}(H)=sp_{0,\f}(H)=sp^B (H)$.

(iv) If  $ W(\r,X) = \f'_{ar}(\r,X)$ and $0$ is an $(\A,V)$-regular point
for $ H$,
  then there is $\delta >0$ and $\psi\in\f(\r)$ such that $\widehat{\psi}\in \h(\r)$, $\widehat{\psi} =1$ on $[-\delta,\delta]$ and $(H*\psi)|\,\jj \in \A$.
\endproclaim

\demo{Proof} (i) If $\jj=\r$ there is nothing to prove so take $\jj=\r_+$. For $\va\in V(\r)$ we have $H * \va  =
\frak{F}*\va  + (H-\frak{F})*\va  $. By Proposition 1.3(ii), $((H-\frak{F}) *\va)|\r_+\in C_0 (\r_+,X)$, so by Proposition 2.1 (i) it follows that $(H * \va) |\r_+ \in \A$ if and only if $(\frak{F} * \va) |\r_+ \in \A$.

(ii)Again we need only  consider the
   case $\jj=\r_+$.  For  $f\in L^1(\r)$ we have $(H*f) |\,\r_+=(\frak{F}*f) |\,\r_+ + \xi$ where $ \xi=((H-\frak{F})* f)|\, \r_+ $.  By (2.3), it follows that  $(\frak{F}*f) |\,\r_ + \in \A$ and by Propositions 1.3 (ii), 2.1 (i) we deduce that    $ \xi\in C_0 (\r_+,X)\st \A$. Hence   $(H*f) |\,\r_+\in \A$.

 (iii) By part (i) we have $sp_{\A,\f} (H)= sp_{\A,\f} (F)$ and $sp_{\A,L^1 (\r)} (H)=sp_{\A,L^1 (\r)} (F)$.
  Moreover, it is clear that $sp_{\A,L^1 (\r)} (F)\st sp_{\A,\f} (F)$. Conversely, we prove  that a point $\om_0\in \r$ is  $(\A,\f)$-regular for $
F$ if  there is $h_0\in L^1(\r)$ such
that $\widehat{h_0}(\om_0)\not =0$ and $(\frak{F}*h_0)|\,\jj \in \A$. Choose $\delta >0$ such that
 $\widehat{h_0}\not = 0$ on $
[\om_0-\delta,\om_0+\delta]$ and by   Lemma 1.1,
  $k_0\in L^1(\r)$ such that $\widehat{k_0}\cdot \widehat{h_0}=1
$ on $ [\om_0-\delta,\om_0+\delta]$. Let $\va\in \f (\r)$, $\widehat {\va}(\om_0)\not =0$ and
supp $\widehat {\va} \st [\om_0-\delta,\om_0+\delta]$.
 By (1.11)  we have $ \frak{F}*\va =
\frak{F}*(h_0*k_0*\va)=(\frak{F}*h_0)*(k_0*\va)$. So, $ (\frak{F}*\va)|\jj \in \A$ by  Proposition 1.3 (i) and part (ii). The second
part follows by taking $\A=\{0\}$.

(iv)  By (i), $0$ is  is an $(\A,\f)$-regular point for $F$, so
 there is $\delta > 0$ and  $\va\in \f(\r)$ such that  $\widehat{\va}\not =0$ on $[-\delta,\delta]$ and $(\frak{F}*\va)|\,\jj \in \A$. If $\jj=\r$, then $H=F=\frak{F}$ and so $(H*\va) \in \A$. If $\jj=\r_+$, then  $H*\va= \frak{F}*\va+ (H- \frak{F})*\va$. So $(H*\va)|\r_+ \in \A$ by Propositions 1.3(ii) and 2.1 (i).
  By Lemma 1.1, there is $g\in \f (\r)$ such that $\widehat{g}\in \h(\r)$ and $\widehat {\psi}= \widehat{\va*g} =1$ on $[-\delta,\delta]$. Obviously $\psi\in \f(\r)$ and $\widehat{\psi}\in \h(\r)$. By (1.11) we have $H*\psi = (H*\va) *g$ and so $(H*\psi)|\,\jj \in \A$ by  part (ii).
$\P$
\enddemo

 \proclaim{Proposition 2.3}  Let  $\A \st L^{\infty} (\jj,{X})$ be a closed subspace  satisfying (2.3).

 (i)  Let   $W,V$ satisfy  (1.8). If $F\in
W(\jj,X)$ and  $g\in V (\r)$ or $g\in L_c^{\infty} (\r)$,  then

(2.5) \qquad $sp_{\A,V} (\frak{F}*g) \st sp_{\A,V} (F) \,\cap $ supp $ \widehat
{g} $  \,  and
  $sp_{\A,V} (F)= \cup_{h > 0} sp_{\A,V} (M_h F)$.

(ii) If   $F\in \f'_{ar}(\jj,X)$, $t\in \r$ and  $0\not = c\in \cc$, then
$sp_{\A,\f} (c\frak {F}_t) = sp_{\A,\f} (F)$.

(iii) If  $F, H\in \f'_{ar}(\jj,X)$, then $sp_{\A,\f} (F+H) \st sp_{\A,\f} (F) \cup sp_{\A,\f} (H) $.

(iv) If  $F\in
\f'_{ar}(\jj,X)$ and $\gamma_{\la} \A \st \A$ for all $\la\in \r$,  then $sp_{\A,\f} (\gamma_{\la}F)=\la +sp_{\A,\f} (F)$.
\endproclaim

\demo{Proof}
(i)  Assume
$\om\not \in sp_{\A,V} (F)$. Then there is $\va\in V (\r)$ with
$\widehat{\va} (\om)\not =0$ and $(\frak{F}*\va)|\jj \in \A$. By (1.9) and Proposition 1.3 (i), $\frak{F}*\va \in BC(\r,X)$.  By (1.11),
 we have $(\frak{F}*g)*\va= (\frak{F}*\va)*g= \frak{F}*(\va*g)$. So, by
Proposition 2.2 (ii), we get  $((\frak{F}*\va)*g)|\,\jj\in \A$
proving $\om\not \in sp_{\A,V} (\frak{F}*g)$. On the other hand if $\om\not
\in $ supp $ \widehat {g}$, then there is $\va\in V (\r)$ with
$\widehat{\va} (\om)\not =0$ and  $\va *g=0$. So, $\om\not \in
sp_{\A,V} (\frak{F}*g)$. For the case $\A=\{0\}$ see also [26, Proposition
0.6 (i)].

To prove  the second  part of (2.5) we note that $M_h F= (\frak{F}*s_h) |\jj$,  $\frak{F}*s_h$  satisfies (2.4) if $\jj=\r_+$ (see  (1.7)) and $g=s_h \in L_c^{\infty} (\r)$  for each $h >0$. Hence $sp_{\A,V} (\frak{F}*s_h)\st sp_{\A,V} (F)$. By Proposition 2.2 (i), we have  $sp_{\A,V} (M_h F)=sp_{\A,V} (F*s_h)$. It follows that
 $ \cup_{h > 0}sp_{\A,V} (M_h F)\st sp_{\A,V} (F)$. Now, let $\om\in sp_{\A,V} (F)$. There is $h >0$ such that $\widehat{s_h} (\om)\not = 0$. Assume that  $\om \not\in sp_{\A,V} (M_h F)= sp_{\A,V} (\frak{F}* s_h)$. By Proposition 2.2 (iv), there is  $\psi\in V(\r)$ such that $\widehat{\psi} (\om)\not = 0$ and $((\frak{F}*s_h)*\psi)|\jj \in \A$. By (1.11), $(\frak{F}*s_h)*\psi=\frak{F}*(s_h*\psi)$. It follows that $(\frak{F}*(s_h*\psi))|\jj \in \A$. Since $s_h*\psi\in V(\r)$ and $\widehat {s_h*\psi} (\om)\not =0$, we conclude that $\om \not\in sp_{\A,V} (F)$, a contradiction which shows $\om \in sp_{\A,V} (M_h F)$.
 This  proves  $sp_{\A,V} (F)\st \cup_{h > 0}sp_{\A,V} (M_h F)$.

 The proofs of (ii), (iii), (iv) are similar to the case $\A=\{0\}$ ([26, Proposition 0.4]). \P
\enddemo

We recall (see [7, p. 118], [8, p. 1007],  [12],  [29]) that a
function $F\in L^1_{loc} (\jj,X)$ is called $ergodic$  if there is
a constant $m(F)\in X$ such that

\qquad \qquad \qquad sup$_{t\in \jj}|| \frac{1}{T}\int_0^T
F(t+s)\,ds -m(F)||\to 0 $ as $T\to \infty$.

\noindent The limit $m(F)$ is called the $mean$ of $F$. The set of
all such ergodic functions will be denoted by $\E(\jj,X)$. We set $\E_{0}(\jj,X)=\{F\in \E(\jj,X): m(F)=0\} $,
$\E_b(\jj,X) =\E(\jj,X)\cap L^{\infty} (\jj,X)$,
$\E_{b,0}(\jj,X)=\{F\in \E_b(\jj,X): m(F)=0\} $, $\E_{ub}(\jj,X)
=\E(\jj,X)\cap BUC(\jj,X)$ and $\E_{u,0}(\jj,X)
=\E_{ub}(\jj,X)\cap \E_{b,0} (\jj,X)$.

 If $F\in L^1_{loc} (\jj,X)$ and
  $\gamma_{\om}F\in \E (\jj,X)$ for some $\om\in\r$, then

 (2.6)\qquad $\gamma_{\om}M_h\, F\in \E (\jj,X)$ \,\, and \,\,
  $M_h\, \gamma_{\om} F\in \E_b (\jj,X)$\,\, for all $h >0$.

\noindent Moreover, if $F\in L^{\infty} (\jj,X)$ and
$\gamma_{\om}F\in \E_b (\jj,X)$ for some $\om\in\r$, then

(2.7) \qquad $\gamma_{\om}(\frak{F}*g)|\,\jj\in \E_{ub}(\jj,X)$ \,\,\, for all
$g\in L^1 (\r)$.

\noindent To prove (2.6), note that

 $M_T \gamma_{\om} M_h  F= \gamma_{\om}M_h
\gamma_{-\om}M_T \gamma_{\om} F$ and $M_T  M_h \gamma_{\om} F= M_h
M_T \gamma_{\om} F$.

\noindent It follows that $\gamma_{\om}M_h F$, $M_h
\gamma_{\om} F\in \E (\jj,X)$ for all $h >0$. By [8, (2.4)], $M_h
\gamma_{\om} F\in C_b (\jj,X)$ and so $M_h \gamma_{\om} F\in \E_b
(\jj,X)$. For (2.7) note that  if $F\in L^{\infty} (\jj,X)$, then
 $M_h F=(\frak{F}*s_h)|\jj$ (see  (1.7))  is  bounded and
uniformly continuous.  So, $\gamma_{\om}M_hF \in \E_{ub} (\jj,X)$
by (2.6). It follows that $ \gamma_{\om}(\frak{F}*g)|\jj \in \E_{ub} (\jj,X)$
for any step function $g$. Since step functions are dense in
$L^1(\r)$, (2.7) follows.

 Also, we note that

 (2.8)\qquad   $ \E_u(\jj,X): =UC(\jj,X)\cap
 \E(\jj,X)=\E_{ub}(\jj,X)$.

\noindent This follows by Proposition 1.2 (ii) using  (2.6) (see also [8, Proposition 2.9]).

Next we recall the definition of the class of slowly
oscillating functions

\qquad $SO(\jj,X)= UC(\jj,X)+ L^1_{loc,0}(\jj,X)$,

\noindent where (see [17, Lemma 1.6], [4, Proposition 4.2.2] for the case $\jj=\r_+$)

\qquad $L^1_{loc,0}(\jj,X)=\{F\in L^1_{loc}(\jj,X): \lim_{|t|\to
\infty, t\in \jj}F(t)=0 \}$

\noindent
It follows that if $ F\in  L^1_{loc,0}(\jj,X)$ and  $\psi\in\f(\r)$, then

 (2.9)\qquad $  F\in \E_0(\jj,X)$, $ \frak{F}\in  L^1_{loc,0}(\r,X)$,

  (2.10)\qquad  $M_h F\in C_0(\jj,X)$  for all $h >0$ and $(\frak{F}*\psi)|\in C_0(\r,X)$.
\proclaim{Lemma 2.4}  If $F\in L^{\infty}(\r,X)$  and  $0\not\in sp_{0,\f} (F)$, then $F\in \E_{b,0} (\r,X)$. If $F\in SO(\r,X)$ and $0\not\in sp_{0,\f} (F)$, then
$F\in \E_{u,0}(\r,X)+ L^1_{loc,0}(\r,X)\st \E_0(\r,X)$.
\endproclaim

\demo{Proof} If $F\in L^{\infty}(\r,X)$ and $0\not\in sp_{0,\f} (F)$, then by Proposition 2.2 (iii) $0\not \in sp^B (F)$.  By [11, Corollary 4.4], $PF\in BUC(\r,X)$, and hence  $F\in \E_{b,0} (\r,X)$. If $F\in SO(\r,X)$,
let $F=\Phi+\xi$ with $\Phi\in  UC(\r,X)$ and $\xi\in L^1_{loc,0}(\r,X)$. By (2.5), we have $sp_{0,\f} (M_h F) \st sp_{0,\f} (F)$ and so $0\not\in sp_{0,\f} (M_h F)$ for all $h > 0$. Since $M_h F\in UC(\r,X) $, we get $M_h F\in BUC(\r,X) $ for all $h  > 0$  by [11, Theorem 4.2]. Since $M_h \xi\in C_0 (\r,X)$ by (2.10), we get $M_h \Phi\in BUC(\r,X)$ for all $h > 0$. This implies $\Phi=\lim_{h\to 0} M_h \Phi\in BUC(\r,X)$ by Proposition 1.2 (ii). Choose $\delta > 0$ and $\va\in \f (\r)$ such that $\widehat {\va}=1$ on $[-\delta,\delta]$.  By (2.9), (2.10) it follows  that $\eta = (\xi- \xi*\va)\in \E_0 (\r,X)$. Set  $\Psi=\Phi +\xi*\va = F- \eta$. Since   $0\not\in sp_{0,\f} (\eta) \cup sp_{0,\f} (\Phi)$, we get  $0\not\in sp_{0,\f} (\Psi)$.  Since by (2.10) $\Psi\in BUC(\r,X)$, we conclude  that $P\Psi\in BUC(\r,X)$, by [11, Corollary 4.4]. This implies that $\Psi\in \E_{u,0} (\r,X)$ and proves that $F=\Psi+\eta \in \E_0 (\r,X)$. \P
\enddemo

We are now ready to state and prove our main results.

\proclaim{Theorem 2.5} Let  $\A \st L^{\infty} (\jj,{X})$ be a closed subspace  satisfying (2.3) and $
\gamma_{\la}\A\st \tilde{\E}\in \{ \E(\jj,X), \E_{0}(\jj,X) \}$ for all $\la\in \r$. Let $F\in \f'_{ar} (\jj,X)$ and $\om\not\in sp_{\A,\f}(F)$.

(i) If $F\in L^{\infty}(\jj,X)$, then  $\gamma_{-\om} F\in \tilde{\E}$.

(ii) If $\gamma_{-\om}F\in SO(\jj,X)$, then $\gamma_{-\om}F \in (\E_{ub}(\jj,X)+ L^1_{loc,0}(\jj,X))\cap\tilde{ \E} $. If also
 $\A= C_0(\jj,X)$, then $\gamma_{-\om} F\in \E_{0}
(\jj,X)$ and if $\gamma_{-\om}F\in UC(\jj,X)$, then $\gamma_{-\om} F\in \E_{u,0}
(\jj,X)$.
\endproclaim

\demo{Proof}  Replacing   $\gamma_{-\om}F$ by $F$, we may assume
$\om=0$ and $0\not\in sp_{\A,\f}(F)$. So, by Proposition 2.2 (iii)    there is  $\delta >0$ and $\va\in
\f(\r)$ such that
 supp $\, \widehat {\va}$ is compact,
 $\widehat {\va} =1$ in a neighbourhood of $0$ and
$(\frak{F}*\va)|\,\jj \in \A\st \tilde{\E} $.
  Set $G=\frak{F}-\frak{F}*\va $.

 (i)  If $F\in L^{\infty}(\jj,X)$, then $G\in L^{\infty}(\r,X)$ and  $0\not\in sp_{0} (G)=
sp_{0,\f} (G)$ by Proposition 2.2 (iii)( see also [9, (3.3), (3.11)]). By Lemma 2.4, we get $G\in \E_{b,0}(\r,X)$. It follows that $F = [G+
{F}* \va]\,|\jj \in \tilde{\E} $.

(ii) If $F\in SO(\jj,X)$, then by (2.10) $G\in SO(\r,X)$  and $0\not\in
sp_{0,\f} (G)$.  By Lemma 2.4, it follows that $G\in \E_{u,0}(\r,X)+L^1_{loc,0}(\r,X)\st \E_0(\r,X)$.
 This implies $F = [G+
{F}* \va]\,|\jj \in \tilde{ \E} $. Obviously if $\A=
C_0(\jj,X)$, then $F\in \E_{0} (\jj,X)$ and if $F\in UC(\jj,X)$, then $ F\in \E_{u,0}
(\jj,X)$. \P
\enddemo

\proclaim{Theorem 2.6} Assume that $F \in \f'_{ar}
(\jj,X)$,  $sp_{C_0(\jj,X),\f}(F)$ is countable and  $\gamma_{-\om} F $
 $\in
\E(\jj,X)$ for all $\om\in sp_{C_0(\jj,X),\f}(F)$.

 (i) If $F \in UC(\jj,X)$, then $F\in AAP(\jj,X)$. If also $sp_{C_0(\jj,X),\f}(F)$
 $=\emptyset$, then $F\in
C_0(\jj,X)$.

 (ii) If $F \in SO(\jj,X)$, then $F\in AP(\jj,X)\oplus
L^1_{loc,0}(\jj,X)$.   If also $sp_{C_0(\jj,X),\f}(F)$

\noindent $=\emptyset$, then $F\in
L^1_{loc,0}(\jj,X)$.

 (iii) If $F=H|\,\jj$ where $H \in L^{\infty}(\r,X)$ and if  $f\in L^1(\r)$, then $(H*f) |\jj\in AAP(\jj,X)$.

(iv) If $sp_{C_0(\jj,X),\f}(F)\not = \emptyset$ and if $\psi\in \f(\r)$, then $(\frak{F}*\psi) |\jj\in
AAP(\jj,X)$.

(v) If  $sp_{C_0(\jj,X),\f}(F)=\emptyset$ and if $\psi\in \f(\r) $ with
$\widehat{\psi}\in \h(\r)$,
then $\frak{F}*\psi \in C_0(\r,X)$.

(vi) If $sp_{C_0(\jj,X),\f}(F) = \emptyset$  and either $F\in \E
(\jj,X)$ or more generally $M_h F\in BC(\jj,X)$ for all $h> 0$ and if  $\psi\in \f(\r) $,
then $\frak{F}*\psi \in C_0(\r,X)$.
\endproclaim
\demo{Proof} (i) First,  we note that $\A=C_0(\jj,X)$ satisfies the assumptions of Theorem 2.5 with $\E =\E_0 (\jj,X)$.  If $0\not\in
sp_{C_0 (\jj,X),\f}(F)$, then by Theorem 2.5, $F \in \E_{u,0}(\jj,X)$. If $0\in sp_{C_0 (\jj,X),\f}(F)$, then from $F\in UC(\jj,X)$
 and (2.8) we get $F \in \E_{ub}(\jj,X)$.
   Let $\tilde{F}\in
BUC(\r,X)$ be an extension of $F$. Since $C_0 (\jj,X)\st
AAP(\jj,X)$ we get $sp_{AAP(\jj,X)}(\tilde{F})\st sp_{C_0
(\jj,X)}(\tilde{F})$. It follows from Proposition 2.2 (iii) that $sp_{AAP(\jj,X)}(\tilde{F})$ is
countable.  By [5, Theorem 4.2.6]  $F =\tilde{F}|\jj \in
AAP(\jj,X)$. If $sp_{C_0(\jj,X),\f}(F)$
 $=\emptyset$, then $\gamma_{\la}F \in \E_{u,0}(\jj,X)$ for all $\la\in \r$. This implies $F\in
C_0(\jj,X)$.

(ii)   Let  $F= u+\xi$, where $u\in UC(\jj,X)$, $\xi\in L^1_{loc,0}(\jj,X)$. We note that  $sp_{C_0 (\jj,X),\f}(\xi)=\emptyset$ by (2.10) and $\gamma_{\la}\xi \in \E(\jj,X)$ for all $\la\in \r$, by (2.9). Also, we have  $sp_{C_0 (\jj,X),\f}(M_h
F)$ is countable by (2.5).   By (2.6), we get $\gamma_{-\om}M_h F\in \E(\jj,X)$ for
all $\om \in sp_{C_0 (\jj,X),\f}(F)$. It follows that $sp_{C_0(\jj,X),\f}(M_h u)$ is countable  and  $\gamma_{-\om}M_h u\in \E(\jj,X)$ for
all $\om \in sp_{C_0 (\jj,X),\f}(F)$.
 So, by part (i), we conclude that
$M_h u \in AAP(\jj,X))$ for  all $h >0$. By Proposition 1.2 (ii), $ u =\lim _{h\to 0} M_h u\in AAP(\jj,X)) $. It follows that $F\in AP(\jj,X)) \oplus L^1_{loc,0}(\jj,X)$. If $sp_{C_0(\jj,X),\f}(F)$
 $=\emptyset$, then $\gamma_{\la}F \in \E_{u,0}(\jj,X)$ for all $\la\in \r$. This implies that $F\in
L^1_{loc,0}(\jj,X)$.

(iii) Let $f\in L^1(\r)$. Then
$\frak{F}*f \in BUC(\r,X)$. By (2.5), we deduce that
$sp_{C_0 (\jj,X),\f}(\frak{F}*f)$ is countable. By (2.7) we find that
$\gamma_{-\om}(\frak{F}*f)|\,\jj \in \E_{ub}(\jj,X)$ for all $\om\in
sp_{C_0 (\jj,X),\f}(\frak{F}*f)$. It follows  that $(\frak{F}*f)|\,\jj \in  AAP(\jj,X))$, by part (i). By Proposition 1.3 (i), we have $ ((H-\frak{F})*f)|\,\jj\in C_0 (\jj,X)$. Hence  $(H*f)|\,\jj \in  AAP(\jj,X))$.

(iv)  Without loss of generality we may assume $0 \in
sp_{C_0 (\jj,X),\f}(F)$. Then $ F\in \E(\jj,X)$ and so by (2.6), $(\frak{F}*s_h)|\,\jj= M_h F\in
\E_b(\jj,X)$  and
$\gamma_{-\om} M_h F\in \E_b(\jj,X)$ for all $h>0$ and $\om \in
sp_{C_0 (\jj,X),\f}(F)$. By (2.5),  $sp_{C_0 (\jj,X),\f} (M_h F)\st
sp_{C_0 (\jj,X),\f}(F)$. Therefore, by part (iii), $((\frak{F}*s_h) *g)|\jj \in
AAP(\jj,X)$ for all $g\in L^1(\r)$. Take $\psi\in \f(\r)$. It
follows that $M_h(\frak{F} *\psi) |\jj=((\frak{F}*s_h) *\psi) |\jj\in AAP(\jj,X)$ and also $(\Delta_h (\frak{F}*\psi))|\jj= (\frak{F}*\Delta_h
\psi)|\jj = (\frak{F}* hM_h \psi')|\jj =(\frak{F}*(s_h *\psi')) |\jj= ((\frak{F}*s_h)*\psi')) |\jj\in AAP(\jj,X)$.  By
[8, Proposition 1.4], one gets $(\frak{F}*\psi)|\,\jj$ is uniformly continuous.
This implies $(\frak{F}*\psi)|\jj= \lim_{h\searrow  0} M_h(\frak{F}*\psi)|\jj\in AAP(\jj,X)$,  by Proposition 1.2 (ii).

(v)  Let $\om \in K=$ supp $\widehat{\psi}$. Since $C_0 (\jj,X)$ satisfies (2.1), (2.3), by Proposition 2.2 (iv), there is
$f^{\om}\in \f(\r)$ such that $\widehat{f^{\om}}$ has compact
support, $\widehat{f^{\om}}=1$ on an open neighbourhood $V(\om)$
of $\om$ and $(\frak{F} *f^{\om})|\,\jj \in C_0(\jj,X)$. Take
$k^{\om}=f^{\om}*g^{\om}$, where $g^{\om} (t)=\overline{f^{\om}
(-t)}$. By  (1.11) and Proposition 2.2(ii),  we conclude that $(\frak{F} *k^{\om})|\,\jj \in C_0(\jj,X)$.
Consider the open covering $\{V(\om): \om\in K\} $. By
compactness, there is a
 finite sub-covering $\{V(\om_1),\cdots, V(\om_n)\}$ of $K$.
 One has $k=\sum_{i=1}^n k^{\om_i}\in \f(\r)$, supp $\widehat{k}$ is compact,
 $\widehat{k}\ge 1$ on $K$ and $(\frak{F}*k) |\jj\in
C_0(\jj,X)$. By Lemma 1.1, there is $h\in \f (\r)$ such that
$\widehat{h}\cdot \widehat{k}=1$ on $K$. Again by (2.3) and Proposition 2.2 (ii), it follows that
 $(\frak{F}*\psi) |\,\jj= ((\frak{F}*k)*h*\psi) |\jj\in C_0(\jj,X)$. By Proposition 1.3 (i), $\frak{F}*\psi \in C_0(\r,X)$.

(vi)  As in part (iv) we conclude that $(\frak{F}*\psi)|\jj $ is uniformly continuous. It follows that $\frak{F}*\psi \in C_0 (\r,X)$ by (2.5), part (i) and Proposition 1.3 (i).     \P
\enddemo

 \proclaim{Proposition 2.7} Assume $ F\in L^1_{loc}(\jj,X)$ and  $sp_{C_0 (\jj,X),\h} (F)=\emptyset$. If $(\frak{F}*\psi)|\jj$ is  uniformly continuous   for some
$\psi\in \h (\r)$,  then
 $\frak{F}*\psi\in C_0 (\r,X)$.
\endproclaim
\demo{Proof}  Let  $\om\in \r$. There is $\va\in \h(\r)$ such that $\widehat {\va} (\om)\not =0$ and $(\frak{F}*\va)|\jj \in
C_0(\jj,X)$. By Proposition 1.3 (i), we get $\frak{F}*\va\in
C_0(\r,X)$.
   By (1.11)
  and Proposition 2.1 (ii), we have
 $(\frak{F}*\psi)* \va= (\frak{F}*\va)*\psi \in
C_0(\r,X)$. By (2.5),   $sp_{C_0 (\r,X),\f}
(\frak{F}*\psi) \st sp_{C_0 (\r,X),\h}
(\frak{F}*\psi)  \st sp_{C_0 (\r,X),\h}
({F})=\emptyset$.  The result follows from Theorem 2.6 (i).
 \P
\enddemo

The following example shows that the assumption of uniform continuity is essential in Proposition 2.7.
\proclaim{Example 2.8} If $F (t)=e^t$ for  $t\in \r $, then  $sp_{C_0 (\jj,X),\h} (F)=\emptyset$   but $(F*\psi)|\jj$ is unbounded for  each  $\psi\in \h (\r)$ with  $
\int_{-\infty}^{\infty} e^{-s}\psi(s)\, ds \not = 0$.
\endproclaim
\demo{Proof}
 For any $\om\in \r$, choose $a >0$ such that
 $\cos\,
\om\,t$ does not change sign on $[0,a]$. Take $\va\in \h(\r)$ such
that $\va
>0 $ on $ (0,a)$ and supp $\va \st [0,a]$. Let $f(t) =\va(t)$ for
$t\ge 0$, $f(t) = - e^{2t}\va(-t)$ for $t < 0$.
  It follows that $f\in \h(\r)$, $F*f =0$ and $\widehat{f}(\om)\not =0$. This means $sp_{C_0 (\jj,X),\h} (F)=\emptyset$.
  Moreover, for $\psi\in \h(\r)$  we have $F*\psi (t)= c e^t$, where $c=
\int_{-\infty}^{\infty} e^{-s}\psi(s)\, ds$.  So,
$(F*\psi)|\jj$ is unbounded if $c\not = 0$. \P
\enddemo

In the following example we calculate the reduced spectra
 of some functions  whose Fourier transforms  may not be regular distributions.

 \proclaim{Example 2.9}  (i) If $ F\in L^p(\jj,X)$ for some $1\le p
< \infty$, then $M_h F\in C_0(\jj,X)$ for all $h>0$ and
$sp_{C_0 (\r,X),V} (F)=\emptyset$, where $V\in \{\h (\r),\f(\r)\}$.

 (ii) Let $F\in
 \E_{ub}(\jj,X)$ and either $F'\in
L^p(\jj,X)$ for some $1\le p < \infty$ or  more generally  $F' \in L^1_{loc}(\jj,X)$ with $M_h
F'\in C_0(\jj,X)$ for all $h >0$. Then $F\in X\oplus C_0 (\jj,X)$ and  $sp_{C_0 (\jj,X)}(F)\st
\{0\}$.
\endproclaim
\demo{Proof} (i)   By H\"{o}lder's inequality,
 $||M_h F (t)||=(1/h )||\int_0^h F (t+s)\, ds||\le
h^{-1/p}$

 \noindent $(\int_0^h ||F(t+s)||^p\,ds)^{1/p})$, so
   $ M_h F \in C_0 (\jj,X)$
for all  $h >0$.  By (1.7) and Proposition 1.3 (i), we get $  \frak{F}*s_h \in C_0 (\r,X)$
for all  $h >0$. So, $sp_{C_0 (\r,X),V} (\frak{F}*s_h)=sp_{C_0 (\jj,X),V} (M_h F)=\emptyset$ for all $h
>0$. Hence
 $sp_{C_0 (\r,X),V} (F)=\emptyset$ by (2.5).

(ii) By part (i)  we have $h M_hF' (\cdot)= F (\cdot+h)-F(\cdot)\, \in C_0
(\jj,X)$ for all $h>0$.  Let $\widetilde {F}\in BUC(\r,X)$ be given by  $\widetilde {F} =F$  on $\jj$ and $\widetilde {F} (t) =F (0)$ on $\r\setminus \jj$.  It follows that $\Delta_s \widetilde {F} \in C_0
(\r,X)$ for all $s\in\r$.  By [5, Theorems
4.2.2, Corollary 4.2.3], we conclude that $F= \widetilde{F}|\jj\in X\oplus C_0 (\jj,X)$. This implies  $sp_{C_0 (\jj,X)}(F)\st
\{0\}$.  $\P$
\enddemo
The following result  is due to Chill
[17, Proposition 2.1]. The proof below is direct and shorter. It follows in particular that the assumptions of Proposition 2.7 are  implied by the assumptions of Example 2.10.

\proclaim{Example 2.10}  If $F \in L^1_{loc}
 (\r,X)$ and   if  the Fourier transform $\widehat{F}$ defined by (1.3) belongs to  $ \f'_{ar}(\r,X)$ and if $\psi\in \h(\r)$, then $F*\psi \in C_0 (\r,X)$ and so $sp_{C_0(\r,X),\h} (F)=\emptyset$.
\endproclaim
\demo{Proof} By (1.1), we have $G= 1/(2\pi) \widehat{F}\widehat {\psi}\in L^1(\r,X)$. By (1.3)

$F*\psi(t)= \int_{-\infty}^{\infty} F (s) \psi(t-s)\, ds= <F,(\check{\psi})_{-t}>=
(1/2\pi)\int_{-\infty}^{\infty} \widehat{F}(\eta) e^{it\eta}\widehat {\psi}(\eta)\, d\eta$.

\noindent This means that  $F*\psi (t)= \widehat {G}(-t)$.  By the Riemann-Lebesgue lemma,  $F*\psi\in C_0(\r,X)$. This implies  $sp_{C_0(\r,X),\h} (F)=\emptyset$. \P
\enddemo

 \head{\S 3. Properties of the weak Laplace spectra}\endhead

In this section we  establish some new properties of the (weak)
Laplace spectrum for  regular tempered distributions and show that
they are similar to those of the Carleman  spectrum (see [26,
Proposition 0.6]). We use the functions $e_a$ for $a \ge 0$ defined on $\r$ or $\r_+$ by $e_a (t)= e^{-at}$.

If $F\in  \f'_{ar}(\r_+,X)$, then
  $ e_a F\in L^1(\r_+,X)$ for all
$a> 0$
 and so the $Laplace$ $ transform$ $\Cal {L}F$  may be defined by

(3.1) \qquad $\Cal {L}F(\la)=\int_0^{\infty}\, e^{-\la \,t}F (t)\,
dt$ \,\,\,\, for\,\,\,\, $\la \in \cc_+ $.

\noindent For a function $F\in  \f'_{ar}(\r,X)$  the $Carleman$ $transform$  $\Cal {C}F$ is defined  by

(3.2) \qquad   $\Cal {C} {F} (\la)=
 \cases { \Cal{L^+}F (\la)= \int_0^{\infty}\, e^{-\la \,t}F (t)\,
dt \qquad \,{\text{for\,\,} \la
 \in \cc_+}}\\
  { \Cal{L^-}F(\la) = - \int _0^{\infty}\, e^{\la\, t} F(-t)\, dt\qquad {\text{\,
  for\,\,}
 \la \in \cc_-}.}\endcases$

If $F\in L^1(\r_+,X)$, then  $\Cal {L} F$ has a continuous
extension  to $\cc_+\cup i\r$  given also by the integral in
(3.1). By the Riemann-Lebesgue lemma $\widehat{\frak{F}} =\Cal
{L}F(i\cdot)\in C_0 (\r,X)$.

 If $F\in  \f'_{ar}(\r_+,X)$, then
 $\widehat {\frak{F}}\in  \f'(\r,X)$ and
$\Cal{L} F(a+i\cdot)=\widehat{e_a\frak{F}}  \in \f'_{ar}(\r,X)$ for all $a >0$.
Moreover, for $\va\in \f(\r)$,

(3.3)\qquad  $<\Cal {L} F (a+i\cdot),\va>\,\,= <\widehat{e_a\frak{F}},\va> =$
 $ <{e_a\frak{F}},\widehat {\va}>\,\,\to
<\frak{F},\widehat {\va}>\,\,= <\widehat {\frak{F}},\va>$,

\noindent where the limit
exists as $a\searrow 0$ by the Lebesgue  convergence
theorem. This means that $\lim_{a\searrow 0}\Cal{L} F (a+i\cdot)= \widehat{\frak{F}} $ with respect to the weak dual topology on $ \f'(\r,X)$.

For a  holomorphic function  $\zeta: \Sigma\to X$, where
$\Sigma=\cc_+$ or $\Sigma=\cc\setminus i\,\r$,  the point
$i\,\om\in i\,\r$ is called a $regular$ $point$ for $\zeta$ or
$\zeta$ is called $holomorphic$ at $i\,\om$, if $\zeta$ has an
extension $\overline{\zeta}$ which is holomorphic in a
neighbourhood $V\st \cc$ of $i\,\om$.

 Points $i\,\om$ which are $not$ $regular$ $points$ are
called $singular$ $points$.

  The  $Laplace$ $spectrum$ of a function    $F\in  \f'_{ar}(\r_+,X)$   is defined by

 (3.4) \qquad $sp^{\Cal{L}}
(F) =\{\om\in\r: i\, \om$ is a singular point for $\Cal {L}F \}$.

\noindent  The $Carleman \,\, spectrum$ of a function    $F\in  \f'_{ar}(\r,X)$  is defined by

(3.5) \qquad $sp^{\Cal{C}} (F) =\{\om\in\r: i\, \om$ is a
singular point for $\Cal {C} {F}\}$. See [4, (4.26)].

\noindent The Laplace spectrum is also called  the half-line
spectrum ([4, p. 275]).

 Note that if $ \overline{\Cal{L}} \gamma_{-\om} F $ and $
\overline{\Cal{C}} \gamma_{-\om} F$  are  holomorphic extensions of
$ \Cal{L} \gamma_{-\om} F$ and $\Cal{C} \gamma_{-\om} F$
respectively, which are holomorphic in a neighbourhood of $0$, then

(3.6)\qquad $\lim _{\la\to 0}\Cal{L} \gamma_{-\om }F (\la)=
\overline{\Cal{L}} \gamma_{-\om} F (0)$  if $\om \not\in
sp^{\Cal{L}} (F) $, and

\qquad\qquad\, $\lim _{\la\to 0}\Cal{C}
\gamma_{-\om }F(\la)= \overline{\Cal{C}} \gamma_{-\om} F (0)$ if
$\om \not\in sp^{\Cal{C}} (F) $.

\noindent If  $F\in L^{\infty} (\r_+,X)$    and $
sp^{\Cal{L}}(F)=\emptyset$, then by Zagier's result [31, Analytic
Theorem] we conclude that $\widehat {F}(\om)=\int_0^{\infty}
e^{-i\om \,t} F (t)\,dt$ exists as an improper  integral (and by (3.6)
equals $\overline{\Cal{L}} \gamma_{-\om} F (0)$) for each $\om \in
\r$. Zagier's Analytic Theorem does not hold  for unbounded
functions. Indeed, the Laplace spectrum of $F (t) = t e^{it^2}$ is
empty  (see Example 3.4 below) and  it can be verified that
$\int_0^{\infty}  e^{-i\om \,t} F (t)\,dt$ does not exist as an
improper Riemann integral for any $\om\in \r$.

For a  holomorphic function  $\zeta: \cc_+\to X$,  the point
$i\,\om\in i\,\r$ is called a  $weakly$ $regular$ $point$ for
$\zeta$ if there exist $\e
> 0$ and  $h\in L^1 (\om-\e, \om+\e)$ such that

(3.7)\qquad  $\lim_{a\searrow 0} \int _{-\infty}^{\infty}\zeta
(a+i\,s)\va (s)\, ds= \int _{\,\om-\e}^{\,\om+\e}h (s)\va (s)\,
ds$

\qquad  \qquad \,\,\, for all $\va\in \h(\r) $ with supp$
\,\va \st (\om-\e, \om+\e)$.

 Points $i\,\om$ which are not weakly regular points are
called $weakly$ $singular$ $points$.

 The $weak\,\, Laplace\,\, spectrum$ of $F\in
   \f'_{ar}(\r_+,X)$ is defined by ([4, p. 324])

(3.8) \qquad $sp^{w\Cal{L}} (F) =\{\om\in\r: i\, \om$ is not a
weakly regular point for $\Cal {L} {F}\}$.

\noindent For  $F\in
   \f'_{ar}(\r,X)$, we write  $sp^{w\Cal{L}^+} (F)= sp^{w\Cal{L}} (F|\r_+)$. It follows readily that if $F\in
   \f'_{ar}(\r,X)$ then

(3.9) \qquad $sp^{w\Cal{L}^+} (F)\st sp^{\Cal{L}^+} (F)\st
sp^{\Cal{C}} (F)$; and,
 if $F\in L^1(\r_+,X)$,
  $sp^{w\Cal{L}} (F)=\emptyset$.

In the following $sp^{*}$ denotes $sp^{\Cal{L}^+}$ or $sp^{w\Cal{L}^+}$  or
$sp^{\Cal{C}}$.

\proclaim{Proposition 3.1}  If $F
\in \f'_{ar}(\r,X)$, then

(i) $sp^{*} (F) =sp^{*} (F_s)$ for each $s\in \r$.

(ii) $sp^{*} (F)=\cup_{h>0} sp^{*} (M_h F)$.

(iii)  $sp^{*} (\gamma_{\om}F)= \om + sp^{*} (F)$.

\endproclaim

\demo{Proof}  (i) A simple calculation shows that  for $\la\in
\cc^{\pm}$

(3.10)\qquad $\Cal {L^{\pm}}F_s (\la) = e^{\la \,s}\Cal {L^{\pm}}F
(\la)- e^{\la \,s} \int_0^{s}\, e^{-\la \,t }F (t)\,dt$.

\noindent Note that the second term on the right of (3.10) is entire in $\la$ for each $s\in \r$. It follows that $\Cal {L^+}F$  (respectively $\Cal {C}F$) is
holomorphic at $i\,\om$ if and only if $\Cal {L}^+ F_s $
(respectively $\Cal {C}F_s$) is holomorphic at $i\,\om$. This
proves (i)   for $sp^{\Cal{L}^+}$ and $sp^{\Cal{C}}$.  Now, assume $i\,\om$
is a weakly  regular point for $\Cal{L}^+F$. So there exists  $\e >0$
and $h \in L^1 (\om-\e,\om+\e)$ satisfying $\lim_{a\searrow 0}
\int _{-\infty}^{\infty}\Cal{L}F (a+i\,\eta)\va (\eta)\, d\eta=
\int _{\om-\e}^{\om+\e}h (\eta)\va (\eta)\, d\eta$
 for all $\va\in \h(\r) $ with supp$ \,\va \st
(\om-\e, \om+\e)$. Then by [28, Theorem 6.18, p. 146] (valid also for $X$-valued distributions),
$\lim_{a\searrow 0} \int _{-\infty}^{\infty}\Cal{L}^+F
(a+i\,\eta) e^{(a+i\,\eta) s}\va (\eta) d\eta= \int
_{\om-\e}^{\om+\e}h (\eta)e^{i\eta \,s} \,\va (\eta) d\eta$
 for all $\va\in \h(\r) $ with supp$ \,\va \st \, (\om-\e, \om+\e)$.
 It follows that $i\,\om$ is a weakly regular point for
$\Cal{L}^+F_s$.

(ii) Another  calculation shows that for  $\la \in \cc^{\pm}$

(3.11)\qquad $\Cal {L^{\pm}}M_h F (\la) =$
  $g(\la\, h) \Cal {L^{\pm}}F (\la) - (1/h)\int_0^h (e^{\la \, v} \int_0 ^v e^{-\la \, t}F(t)
 dt)\,dv$,

\noindent where  $g$ is the entire function given by
$g(\la)=\frac{e^{\la}-1}{\la}$ for $\la\not = 0$. Let $i\,\om\in
i\, \r$ be a regular point for $\Cal {L^+}F$ and let $\overline
{\Cal {L^+}}F: V \to X$ be a holomorphic extension of $\Cal
{L^+}F$ to a neighbourhood $V \st \cc$ of $i\,\om$. Then
$\overline {\Cal {L^+}}M_h F (\la) =$
  $g(\la\,h)\overline  {\Cal {L^+}}F (\la) -
  (1/h)\int_0^h (e^{\la \, v} \int_0 ^v e^{-\la \, t}F(t)\,
 dt) \,dv$, $\la \in V$,  is a holomorphic extension of
$\Cal {L^+}M_h F$. So $i\,\om$ is a regular point for $\Cal {L^+}
M_h F$.
 Conversely suppose  $i\om\in i\r$  is a regular point of $\Cal{L}^+ M_h F$ for each $h> 0$. Choose $h_0 >0$ such
that $g(i\om\, h_0)\not =0$. Then $i\,\om$ is
a regular point for $\Cal {L^+} F$. This proves (ii) for
$sp^{\Cal{L}^+}$. The case $sp^{\Cal{C}}$ follows similarly noting
that (3.11) implies $\Cal {C} M_h F (\la) =$
  $g(\la\, h) \Cal {C}F (\la) - (1/h)\int_0^h (e^{\la \, v} \int_0 ^v e^{-\la \, t}F(t)
 dt)\,dv$. The proof for  $sp^{w\Cal{L}^+}$ is
similar to the one in part (i).

(iii) This follows easily from the definitions noting that $\Cal
{L}^+(\gamma_{\om}F) (\la)= \Cal {L}^+ F (\la -i\om)$  and $\Cal
{C}(\gamma_{\om}F) (\la)= \Cal {C}F (\la -i\om)$.
  \P
\enddemo

Proposition 3.1 holds for $\frak{F}$, where $F
\in \f'_{ar}(\r_+,X)$. In this case   $sp^{\Cal{L}^+} \frak{F}= sp^{\Cal{L}} F$ and $sp^{w\Cal{L}^+} \frak{F}= sp^{w\Cal{L}} F$.

The following result was obtained in  [16, Lemma 1.16] in the case $\A= C_0(\r_+,X)$ and $F \in L^{\infty}(\r_+,X)$ since then $sp_{C_0(\r_+,X),\f}(F)=sp_{C_0(\r_+,X)}(F)$.

\proclaim{Proposition 3.2}  If $F \in \f'_{ar}(\r_+,X)$ and $\A\st L^{\infty}(\r_+,X)$ satisfies (2.1) then $ sp_{\A,\f}(F) \st sp_{C_0(\r_+,X),\f}(F) \st sp^{w\Cal{L}}(F)$.
\endproclaim
 \demo{Proof}  By Proposition 2.1 (i),  $C_0(\r_+,X)\st \A$ and so $ sp_{\A,\f}(F) \st sp_{C_0(\r_+,X),\f}(F) $.  Let $\om \not \in  sp^{w\Cal{L}}(F) $. Choose $\e >0$ and
$\va\in\f(\r)$ such that $sp^{w\Cal{L}}(F)\cap
[\om-\e,\om+\e]=\emptyset$, $\widehat{\va}(\om)=1$ and supp $\widehat{\va}\st
[\om-\e,\om+\e]$. By [17, Proposition 1.3], $\frak{F}*\va\in C_0(\r,X)$ and so $\om\not\in sp_{C_0(\r_+,X),\f}(F)$. \P
\enddemo

 \proclaim{Remark 3.3} (i)  In the case $\jj=\r_+$  Theorem 2.5 and Theorem 2.6 remain valid if we replace $sp_{\A,\f} (F)$ and $sp_{C_0(\r_+,X),\f} (F)$ by $sp^{\Cal{L}} (F)$ or $sp^{w\Cal{L}} (F)$.
 Indeed,  note that Theorem 2.5  holds for $\A= C_0 (\r_+,X)$. By   Proposition 3.2 and (3.9), we have $sp_{C_0(\r_+,X),\f} (F)\st sp^{w\Cal{L}} (F) \st sp^{\Cal{L}} (F)$.

(ii) If $F$ in Theorem 2.5  is
not bounded or slowly oscillating, then $F$ is not necessarily ergodic.
For example, if $\frak{g} (t)= e^{it^2}$ and $F= \frak{g}^{(n)}$ for some $n\in \N$, then by
Example 3.4 below  and (3.9), we find $sp^{w\Cal{L}^+} (F)=\emptyset$. By
Proposition 3.2, we get $sp_{C_0 (\r_+,\cc)} (F )=\emptyset$ but
$F|\r_+ $ is neither bounded nor  ergodic  when $n\ge 2$. If $n=1$,
$F$ is ergodic but  not bounded.

 (iii) In view of Proposition 3.2 and (3.9) several tauberian theorems by Ingham  [21] ([4, Theorem 4.9.5]) and their
generalizations in [2], [3], [4, Theorem 4.7.7, Corollary 4.7.10, Theorem 4.9.7, Lemma 4.10.2], [13], [14], [16, Lemma 1.16, p. 25], [17] are consequences of
Theorem 2.5 and Theorem 2.6. Our
proofs are simpler and  different. Replacing Laplace and weak Laplace spectra by reduced spectra we are able to strengthen and unify these previous results.
\endproclaim

In the following  we use our results to calculate some (weak) Laplace spectra.

\proclaim{Example 3.4} Take $\frak{g} (t)=e^{it^2}$ for  $t\in \r$. Then  $sp^{\Cal {C}}(\frak{g}) =\r$ and
$sp^{\Cal {L}^+}(\frak{g})= sp^{\Cal {L}^+}(\frak{g}^{(n)}) =\emptyset$ for any
$n\in \N$. Moreover, $M_h \frak{g} \in C_0
(\r,\cc)$ and $sp^{\Cal {L}^+}(M_h \frak{g}) =\emptyset$ for all  $h>0$.
\endproclaim
\demo{Proof} By Proposition 3.1 (i), (iii), it is readily verified
 that $sp^{\Cal {L}^+}(\frak{g})=sp^{\Cal {L}^+}(\frak{g}_a)=2a+ sp^{\Cal
 {L}^+}(\frak{g})$ for each  $a\in \r$. This implies that either $sp^{\Cal {L}^+}(\frak{g})
 =\emptyset$ or $sp^{\Cal {L}^+}(\frak{g}) =\r$. Similarly either $sp^{\Cal {C}}(\frak{g})
 =\emptyset$ or $sp^{\Cal {C}}(\frak{g}) =\r$.  But $\frak{g}\not =0$ and so by [26, Proposition 0.5 (ii)], $sp^{\Cal {C}}(\frak{g})=\r$.

  Next note that $ y (\la)= \Cal
 {L}^+\frak{g} (\la)$ is a solution of the differential equation $y'(\la) + (\la/2i)y (\la)= 1/2i$ for $\la\in \cc_+$. Solving the equation we find $y (\la)= e^{-\la ^2/4i} (c+ (1/2i) \int_0^{\la} e^{z^2/4i}\,dz$ for some choice of $c\in \cc$. As this last function is entire we conclude that $ sp^{\Cal {L}^+}(\frak{g})=\emptyset $.
 Since  $\int_0^{\infty} e^{i\,t^2}\, dt$ converges as
an improper Riemann integral and $M_h \frak{g}(t)= P\frak{g}
(t+h)-P\frak{g} (t)$ it follows that  $M_h \frak{g}\in C_0 (\r,\cc)$ for each $h>0$. Moreover, by Proposition 3.1(ii),  $sp^{\Cal {L}^+}(M_h \frak{g})
=\emptyset$.
 \P
\enddemo

Finally, we demonstrate  that our results can be used to deduce spectral criteria for bounded solutions of evolution equations  of the form

 (3.12)\qquad  $\frac{d u(t)}{dt}= A u(t) +  \phi (t) $, $u(0) \in X$, $t\in  {\jj}$,

\noindent  where
 $A$ is a closed linear operator   on $X$ and $\phi\in L^{\infty} (\jj, X)$.

\proclaim{Theorem 3.5} Let  $\phi\in  L^{\infty} (\jj, X)$ and $u$ be a bounded mild solution of (3.12). Let $\A $ satisfy (2.1), (2.3), $\gamma_{\la}\A \st \A$ for all $\la\in \r$ and contain all constants.

(i) If $\jj=\r_+$, then $i\, sp^{\Cal {L}} (u)\st (\sigma(A)\cap i\r)\cup i\,sp^{\Cal {L}}(\phi)$.

(ii) If $sp_{\A} (\phi)=\emptyset$, then
  $i\,sp_{\A} (u) \st \sigma(A)\cap i\r$.
\endproclaim

\demo{Proof}  As $u, \phi\in L^{\infty} (\jj,X)$ we get  $M_h u, M_h \phi\in BUC(\jj,X)$  and $v= M_h u$ is a classical solution of $ v'(t)= A v(t) + M_h  \phi (t) $, $v(0) \in D(A)$, $t\in  {\jj}$ for each $h >0$.

(i)  By [4, Proposition 5.6.7, p. 380], we have

$i\,sp^{\Cal {L}} (M_h u)\st (\sigma(A)\cap i\r)\cup i\,sp^{\Cal {L}} (M_h\phi)$ for all $h >0$.

\noindent Taking the union of both sides, we get

 $\cup _{h >0}\, i\,sp^{\Cal {L}} (M_h u)\st (\sigma(A)\cap i\r)\cup (\cup _{h >0} \, i\,sp^{\Cal {L}} (M_h (\phi))$.

\noindent Applying Proposition 3.1 (ii) to both sides, we conclude that

$i\,sp^{\Cal {L}} (u)\st (\sigma(A)\cap i\r)\cup i\,sp^{\Cal {L}} (\phi)$.

 (ii)  Take $h >0$. Since $sp_{\A} (\phi)=\emptyset$, it follows that $sp_{\A} (M_h\phi)=\emptyset$,  by (2.5). Hence $M_h\phi\in \A$ by [5, Theorem 4.2.1].  Using [6, Corollary 3.4 (i)], we conclude that  $i\,sp_{\A} (M_h u) \st \sigma(A)\cap i\r$. Again by (2.5), we conclude that $i\,sp_{\A} (u) \st \sigma(A)\cap i\r$. For further details see [10, Proposition 4.2, Theorem 4.3].
  \P
\enddemo

\smallskip

\smallskip

\indent School of Math. Sci., P.O. Box 28M, Monash University,
 Vic. 3800.

\indent Email "bolis.basit\@monash.edu",\qquad
"alan.pryde\@monash.edu".

\enddocument